\input vanilla.sty
\scaletype{\magstep1}
\scalelinespacing{\magstep1}
\def\bull{\vrule height .9ex width .8ex depth -.1ex}

\title Interpolation of compact operators by \\the methods
of Calder\'on and Gustavsson-Peetre \endtitle

\author M. Cwikel\footnote{Research supported by US-Israel
BSF-grant 87-00244 and by the Technion V.P.R.  Fund - B. and
G. Greenberg Research Fund (Ottawa)} \\ Department of
Mathematics\\ Technion, Israel Institute of Technology, \\
Haifa 32000, Israel\\ {\it and}\\ N.J.
Kalton\footnote{Research supported by NSF-grants DMS-8901636
and DMS-9201357; the author also acknowledges support from
US-Israel BSF-grant 87-00244} \\ Department of Mathematics\\
University of Missouri-Columbia\\ Columbia, Missouri 65211,
USA.  \endauthor

\vskip2truecm

\subheading{Abstract}Let $\bold X=(X_0,X_1)$ and $\bold
Y=(Y_0,Y_1)$ be Banach couples and suppose $T:\bold X\to
\bold Y$ is a linear operator such that $T:X_0\to Y_0$ is
compact.  We consider the question whether the operator
$T:[X_0,X_1]_{\theta}\to [Y_0,Y_1]_{\theta}$ is compact and
show a positive answer under a variety of conditions.  For
example it suffices that $X_0$ be a UMD-space or that $X_0$
is reflexive and there is a Banach space so that
$X_0=[W,X_1]_{\alpha}$ for some $0<\alpha<1.$

\vskip2truecm

\subheading{1.  Introduction}

Let ${\bold X}=(X_{0},X_{1})$ and ${\bold Y}=(Y_{0},Y_{1})$
be Banach couples and let $T$ be a linear operator such that
$T:{\bold X}\rightarrow {\bold Y}$ (meaning, as usual, that
$T:X_{0}+X_{1}\rightarrow Y_{0}+Y_{1}$ and
$T:X_{j}\rightarrow Y_{j}$ boundedly for $j=0,1$).
Interpolation theory supplies us with a variety of {\it
interpolation functors} $F$ for generating {\it
interpolation spaces}, i.e. functors $F$ which when applied
to the couples ${\bold X}$ and ${\bold Y}$ yield spaces
$F({\bold X})$ and $F({\bold Y})$ having the property that
each $T$ as above maps $F({\bold X})$ into $F({\bold Y})$
with bound $$ \| T\| _{F({\bold X})\rightarrow F({\bold
Y})}\le C \max(\| T\| _{X_{0}},\| T\| _{X_{1}}) $$ for some
absolute constant $C$ depending only on the functor F. It
will be convenient here to use the customary notation $\|
T\| _{{\bold X}\rightarrow {\bold Y}}= \max(\| T\|
_{X_{0}},\| T\| _{X_{1}})$ . Further general background
about interpolation theory and Banach couples can be found
e.g. in [1], [3] or [5].

In this paper we shall be concerned with the following
question:  \proclaim{Question 1}Suppose that the operator
$T:{\bold X}\rightarrow {\bold Y}$ also has the property
that $T:X_{0}\rightarrow Y_{0}$ is compact.  Let $F$ be some
interpolation functor.  Does it follow that $T:F({\bold
X})\rightarrow F({\bold Y})$ is compact?  \endproclaim

The first positive answer to a question of this type was
given by Krasnosel'ski\v{i} [17] in 1960 in the special
context of $L_{p}$ spaces.  Subsequently Lions and Peetre
[18] gave the first in a series of partial answers to
Question 1 for the case where $F$ is interpolation functor
of the real method $F(\bold X)=(X_0,X_1)_{\theta,p}$, and in
1992 one of us gave the complete answer for the real method
[12], using results and methods suggested by the work of
Hayakawa [15] and Cobos-Peetre [10].

Question 1 is still open in the case where $F$ is the
functor $F(\bold X)=[X_0,X_1]_{\theta}$ of Calder\'on's
complex method [8].  Among the partial solutions which have
been given to date we mention results of Calder\'on [8]
(Sections 9.6, 10.4), Persson [21], Cwikel [12],
Cobos-K\"uhn-Schonbek [9] and a forthcoming paper of Masty\l o
[19].  In this paper we present some further partial results
for this functor.  We are able to answer Question 1 in the
affirmative in each of the following four cases:

\item{(i)}if $X_0$ has the UMD property, \item{(ii)}if $X_0$
is reflexive and is given by $X_0=[W,X_1]_{\alpha}$ for some
Banach space $W$ and some $\alpha \in (0,1)$, \item{(iii)}if
$Y_0$ is given by $Y_0=[Z,Y_1]_{\alpha}$ for some Banach
space $Z$ and some $\alpha \in (0,1)$, \item{(iv)}if $X_0$
and $X_1$ are both complexified Banach lattices of
measurable functions on a common measure space.

Our result (iv) strengthens Theorem 3.2 of [9] where both
$\bold X$ and $\bold Y$ are required to be such couples of
complexified lattices with some other mild requirements.  We
obtain (iv) as a corollary of the result that $T:\langle
X_0,X_1,\theta \rangle \rightarrow [Y_0,Y_1]_{\theta}$ is
compact for arbitrary Banach couples $\bold X$ and $\bold
Y$.  Here $\langle X_0,X_1,\theta \rangle \
=G_{3,\theta}(\bold X)$ denotes the interpolation space
defined by Gustavsson-Peetre [13] and characterized as an
orbit space by Janson [16].  Masty\l o [19] has obtained an
alternative proof of (iv) as a consequence of other results
of his which answer Question 1 in the cases where $F$ is the
Gustavsson-Peetre functor, or other related functors
introduced by Peetre and by Ovchinnikov.

We now recall the definitions of the main interpolation
functors to be used in this paper:

\vskip1truecm

\noindent {\it 1. Calder\'on's complex method
$[\cdot,\cdot]_{\theta}$}

For each Banach couple $\bold X=(X_0,X_1)$ we let $\Cal
H=\Cal H(\bold X)$ denote the space of all $X_0+X_1$-valued
functions which are analytic on the open annulus $\Omega =
\{z:1<|z|<e\}$ and continuous on the closure of $\Omega $.
This space is normed by $\|f\|_{\Cal H}=\max_{z\in
\Omega}\|f(z)\|_{X_0+X_1}.$

The space $\Cal F=\Cal F(\bold X)$ is defined to be the
subspace of $\Cal H$ which consists of those functions $f$
which are $X_0$-valued and $X_0$-continuous on the circle
$|z|=1$ and $X_1$-valued and $X_1$-continuous on the circle
$|z|=e.$ We define $\|f\|_{\Cal
F}=\max_{j=0,1}(\max_{|z|=e^j}\|f(z)\|_{X_j}).$ For each
$\theta \in [0,1]$ the interpolation space
$[X_0,X_1]_{\theta}$ generated by Calder\'on's complex
method is the set of all elements $x\in X_0+X_1$ of the form
$x = f(e^{\theta})$ where $f\in \Cal F$.  Its norm is given
by $\|x\|_{X_{\theta}}=\inf\{\|f\|_{\Cal F}:f\in \Cal F,
f(e^{\theta})=x\}.$ In fact this definition differs slightly
from the one given in Calder\'on's classical paper [8] where
the unit strip $\{z:0<\Re{z}<1\}$ replaces the annulus
$\Omega$ but, as shown in [11], the two definitions coincide
to within equivalence of norms.

We will sometimes use the notation $X_{\theta} =
[X_0,X_1]_{\theta}$ when there is no danger of confusion.
Actually this could be ambiguous for the (sometimes
forgotten) endpoint values $\theta =j=0,1$ since then
$[X_0,X_1]_j$ is the closure of $X_0 \cap X_1$ in $X_j$.
(see [3] Theorem 4.2.2 p. 91 or [8] Sections 9.3 and 29.3,
pp. 116, 133-4.)

The couple $\bold X$ is said to be {\it regular} if $X_0
\cap X_1$ is dense in $X_0$ and also in $X_1$.  If $\bold X$
is regular then the dual spaces $X^*_0$ and $X^*_1$ also
form a Banach couple.  Calder\'on's duality theorem ([8]
sections 12.1, 32.1) states that for regular couples and
$\theta \in (0,1)$ the dual of $[X_0,X_1]_{\theta}$
coincides with the space $[X^*_0,X^*_1]^{\theta}$ obtained
by applying a variant of Calder\'on's construction to the
couple $(X^*_0,X^*_1)$.  We refer to [8] for the exact
definition of this second Calder\'on method
$[\cdot,\cdot]^{\theta}$.

\vskip1truecm

\noindent {\it 2. Peetre's method $\langle \cdot,\cdot
\rangle_{\theta}$}

For each Banach couple $\bold X$ and each $\theta \in (0,1)$
the space $\langle X_0,X_1 \rangle_{\theta}$ is the set of
all elements $x \in X_0+X_1$ which are sums of the form
$x=\sum_{k\in\bold Z}x_k$ where the elements $x_k \in X_0
\cap X_1$ are such that $\sum_{k\in\bold Z}e^{-\theta k}x_k$
is unconditionally convergent in $X_0$ and $\sum_{k\in\bold
Z}e^{(1-\theta)k}x_k$ is unconditionally convergent in
$X_1$.  $\langle X_0,X_1 \rangle_{\theta}$ is normed by $$
\|x\|_{\langle X_0,X_1 \rangle_{\theta}} = \inf \max_{j=0,1}
\sup \|\sum_{k\in\bold Z}\lambda_k
e^{(j-\theta)k}x_k\|_{X_j} $$ where the supremum is taken
over all complex valued sequences $(\lambda_k)$ with
$|\lambda_k|\le 1$ for all $k$, and the infimum is taken
over all representations as above $x=\sum_{k\in\bold Z}x_k$.
We refer to [20] and [16] for more details.

\vskip1truecm

\noindent {\it 3. Gustavsson-Peetre's method $\langle
\cdot,\cdot,\theta \rangle$}

This is a variant of Peetre's method (see [13] and [16]).
The space $\langle X_0,X_1,\theta \rangle$ is defined like
$\langle X_0,X_1 \rangle_{\theta}$ except that the series
$\sum_{k\in\bold Z}e^{(j-\theta)k}x_k$ need only be {\it
weakly} unconditionally Cauchy in $X_j$.  The norm is
accordingly given by $$ \|x\|_{\langle X_0,X_1,\theta
\rangle} = \inf \max_{j=0,1} \sup \|\sum_{k\in F}\lambda_k
e^{(j-\theta)k}x_k\|_{X_j} $$ where the supremum is over all
$\lambda_k$'s as before and over all finite subsets $F$ of
$\bold Z$.

\vskip1truecm

Finally, we discuss the class of UMD-spaces.
Let $X$ be a Banach space and let $\bold T$ denote the unit
circle with normalized Haar measure $dt/2\pi.$ If $f\in
L_2(\bold T,X)$ we denote its Fourier coefficients $$\hat
f(k)=\int_0^{2\pi} f(e^{it})\frac{dt}{2\pi}.$$ Then the
(formal) Fourier series of $f$ is $f\sim\sum_{k\in Z}\hat
f(k)z^k.$ We recall that $X$ is a UMD-space if the
vector-valued Riesz projection $\Cal R:L_2(\bold T,X)\to
L_2(\bold T,X)$ is bounded where $\Cal Rf\sim\sum_{k\ge
0}\hat f(k)z^k.$ In fact UMD-spaces were introduced by
Burkholder in [6] with a different definition, but the above
characterization follows from results of Burkholder [7] and
Bourgain [4].

It perhaps important to stress that although the condition
of being a UMD-space is fairly stringent many of the
well-known spaces used in analysis are in fact UMD.  The
spaces $L_p$ and the Schatten ideals $\Cal C_p$ for
$1<p<\infty$ are UMD; further examples are reflexive Orlicz
spaces and the Lorentz spaces $L(p,q)$ where $1<p,q<\infty$
(see [14]).  The class of UMD-spaces is closed under
quotients, duals and subspaces.  All UMD-spaces are
superreflexive but the converse is false even for lattices
[4].

\vskip2truecm

\subheading{2.  Some preliminary results}

We will make repeated use of the following simple lemma.

\proclaim{Lemma 1}Let $X$ and $Y$ be Banach spaces and
suppose $T:X\to Y$ is a compact operator.  Suppose $(f_n)$
is a bounded sequence in $L_2(\bold T,X)$.  Let $H$ be the
subspace of all elements $y^*\in Y^*$ which satisfy $$
\lim_{n\to\infty} \int_0^{2\pi}|\langle Tf_n(e^{it}),y^*
\rangle |^2 \frac{dt}{2\pi}=0.$$ Suppose $H$ is
weak$^*$-dense in $Y^*$ (i.e.  $H$ separates the points of
$Y$).  Then $$
\lim_{n\to\infty}\int_0^{2\pi}\|Tf_n(e^{it})\|_Y^2
\frac{dt}{2\pi}=0.$$ \endproclaim

\demo{Proof}Let $(y_m^*)$ be a sequence in $H\cap B_{Y^*}$
such that $(T^*y_m^*)_{m=1}^{\infty}$ is norm dense in
$T^*(H\cap B_{Y^*}).$ Then any bounded sequence $(x_n)$ in
$X$ such that $\lim_{n\to\infty}\langle Tx_n,y_m^*\rangle
=0$ for each $m$ must satisfy $\lim_{n\to\infty }\langle
Tx_n,y^*\rangle =0$ for all $y^*\in H$.  Consequently, by
compactness, $\lim\|Tx_n\|_Y=0$.  From this it follows
easily that for every $\epsilon>0$ there exists a constant
$C=C(\epsilon)$ such that $$\|Tx\|^2_Y \le \epsilon\|x\|^2_X
+ C\sum_{m=1}^{\infty}2^{-m} |\langle Tx,y_m^*\rangle|^2$$
for every $x\in X$.  Now the lemma follows
easily.\bull\enddemo

We shall need the following properties of the complex
interpolation spaces $X_{\theta}$, most of which are well
known.

\proclaim{Lemma 2}(i) For each $0<\theta<1$ there is a
constant $C=C(\theta)$ such that, for all $f\in\Cal F$, $$
\|f(e^{\theta})\|_{X_{\theta}} \le C \left(\int_0^{2\pi}
\|f(e^{it})\|_{X_0}\frac{dt}{2\pi}\right)^{1-\theta}
\left(\int_0^{2\pi}
\|f(e^{1+it})\|_{X_1}\frac{dt}{2\pi}\right)^{\theta}.
\leqno (1) $$ In particular, for all $x\in X_0\cap X_1$, $$
\|x\|_{X_{\theta}} \le
C\|x\|_{X_0}^{1-\theta}\|x\|_{X_1}^{\theta}.  \leqno (2) $$

(ii) For each $0<\theta<1$, $X_0\cap X_1$ is a dense
subspace of $X_{\theta}$.

(iii) Let $X^{\circ}_j$ denote the closed subspace of $X_j$
generated by $X_0 \cap X_1$.  Then, for all $\theta \in
[0,1]$, $$ [X^{\circ}_0,X^{\circ}_1]_{\theta} =
[X_0,X^{\circ}_1]_{\theta} = [X^{\circ}_0,X_1]_{\theta} =
[X_0,X_1]_{\theta}$$.

(iv) (reiteration formulae) $$ [[X_0,X_1]_{\theta
_0},[X_0,X_1]_{\theta _1}]_{\sigma } = [X_0,X_1]_s \leqno
(3)$$ with equivalence of norms, for each $\theta _0$,
$\theta _1$ and $\sigma $ in $[0,1]$, where $s=(1- \sigma
)\theta _0 + \sigma \theta _1$.  Also $$ [[X_0,X_1]_{\theta
_0},X_1]_{\sigma } = [X_0,X_1]_{(1- \sigma )\theta _0 +
\sigma} \leqno (4)$$ and $$ [X_0,[X_0,X_1]_{\theta
_1}]_{\sigma } = [X_0,X_1]_{\sigma \theta _1}.  \leqno (5)$$
\endproclaim

\demo{Proof}Part (i) follows easily from the above-mentioned
equivalence of complex interpolation in the annulus with
complex interpolation in the unit strip, by applying the
estimate (ii) of [8] section 9.4 p. 117 to the function
$F(z)=f(e^z)e^{z^2}$.

For parts (ii) and (iii) we refer to [3] Theorem 4.2.2 p. 91
or [8] Sections 9.3 and 29.3, pp. 116, 133-4.  For part (iv)
the formula $(3)$ is proved in [11] pp. 1005-1006, and also
in [16], Theorem 21, pp. 67-68.  Its variant $(4)$ follows
from $(3)$ if $X_0 \cap X_1$ is dense in $X_1$.  But it can
also be shown in general by slightly modifying Janson's
proof of the reiteration formula ([16], Theorem 21, pp.
67-68):  One of the things to bear in mind for that proof is
that simple estimates with the $K$-functional show that
$[X_0,X_1]_{\theta _0}\cap X_1 \subset [X_0,X_1]_{(1- \sigma
)\theta _0 + \sigma}$.  (Cf.  [11]).  The proof of $(5)$ is
exactly analogous.  \bull \enddemo

For each $f\in\Cal H$ we write $f(z)=\sum_{k\in\bold Z}\hat
f(k)z^k$ and we let $\Cal Rf$ denote the analytic function
on $\Omega$ defined by $\Cal Rf(z)=\sum_{k\ge 0}\hat
f(k)z^k.$ We set $\Cal R_-f=f-\Cal Rf.$ It is easy to see
that $\Cal Rf$ extends to an $X_0+X_1-$valued analytic
function on the open disk $|z|<e$ and similarly $\Cal R_-f$
is analytic on the open set $|z|>1.$ It thus follows that
$\Cal Rf$ extends to an element of $\Cal H$ and that $\|\Cal
Rf\|_{\Cal H}\le C\|f\|_{\Cal H}$ for some absolute constant
$C.$

For each positive integer $N$ and $f\in \Cal H$ we define
$\Cal S_Nf$ by the formula $$\Cal S_N(f) = \sum_{|k|\le
N}\hat f(k)z^k + \sum_{N<|k|\le 2N}(2-\frac{|k|}N)\hat
f(k)z^k.$$ By the uniform $L_1-$boundedness of the de la
Vall\'ee Poussin kernels there exists a constant $C$ such
that $\|\Cal S_Nf\|_{\Cal F}\le C\|f\|_{\Cal F}$ for all
$f\in \Cal F$ and all $N>0$.

Now let $\bold Y=(Y_0,Y_1)$ be another Banach couple and let
$T:{\bold X}\to {\bold Y}$ be a linear operator with the
further property that $T:X_0 \to Y_0$ is compact.  We may
assume that $\|T\|_{\bold X\to\bold Y}\le 1$.  In fact $T$
will be assumed to have these properties throughout the
remainder of this paper.

\proclaim{Lemma 3}(a) The set $\{T\hat f(k):f\in B_{\Cal
F},\ k\in\bold Z\}$ is relatively compact in $Y_0.$\newline
(b) We have $\lim_{k\to \infty}\sup_{f\in B_{\Cal F}}
\|T\hat f(k)\|_{Y_0} =0.$ \newline (c) For each $\delta>0$
there exists an integer $L=L(\delta)$ so that for each $f\in
B_{\Cal F}$ the set $\{k:\|T\hat f(k)\|_{Y_0}>\delta\}$ has
at most $L$ members.\newline (d) For each $0<\theta<1$ we
have $$ \lim_{|k|\to\infty}\sup_{f\in B_{\Cal F}}\|T\hat
f(k) e^{k\theta}\|_{Y_{\theta}}=0.$$ \endproclaim

\demo{Proof} (a) We simply observe that $$ T\hat f(k) =
T(\int_0^{2\pi} f(e^{it})e^{-ikt}\frac{dt}{2\pi}).$$

(b) Since $T:X_0\to Y_0$ is compact, there exists a function
$\eta:[0,\infty)\to [0,\infty)$ with $\lim_{\delta\to
0}\eta(\delta)=\eta(0)=0$ such that if $\|x\|_{X_0}\le 1$
then $\|Tx\|_{Y_0}\le \eta(\|x\|_{X_1})$ whenever
$\|x\|_{X_1}<\infty.$

Now for $f\in B_{\Cal F}$ we have $\|\hat f(k)\|_{X_0}\le 1$
and $\|\hat f(k)e^k\|_{X_1}\le 1.$ Hence $\|T\hat
f(k)\|_{Y_0} \le \eta(e^{-k}).$

(c) Since $T$ is compact we can pick a finite set of
functionals $\{y_1^*,\ldots,y_N^*\}$ in $B_{Y_0^*}$ such
that for $x\in X_0$ we have $$\|Tx\|_{Y_0}\le \max_{1\le
j\le N} |\langle Tx,y_j^*\rangle | +
\frac12\delta\|x\|_{X_0}.$$ Now suppose $f\in B_{\Cal F}$
and that $A=\{k:\ \|T\hat f(k)\|_{Y_0}>\delta\}.$ Then for
each $k\in A$ we have $$ \sum_{j=1}^N|\langle T\hat
f(k),y_j^*\rangle |^2 \ge \frac14\delta^2.$$ Summing over
all $k$ and applying Parseval's identity, we obtain $$
\sum_{j=1}^N\int_0^{2\pi} |\langle Tf(e^{it}),y_j^*\rangle
|^2 \frac{dt}{2\pi} \ge \frac14 \delta^2 |A|.$$ Thus $|A|\le
4N\delta^{-2}$ from which the result follows immediately.

(d) First we observe, using $(2)$, that $$\|T\hat
f(k)\|_{Y_{\theta}} \le C\|T\hat
f(k)\|_{Y_0}^{1-\theta}\|T\hat f(k)\|_{Y_1}^{\theta} \le
C\|T\hat f(k)\|_{Y_0}^{1-\theta} e^{-k\theta} $$ and so we
obviously have from (b) that $$ \lim_{k\to\infty} \sup_{f\in
B_{\Cal F}}\|T\hat f(k)e^{k\theta}\|_{Y_{\theta}}=0.$$

It remains to establish a similar result as $k\to -\infty.$
Suppose then that this is false.  Then we can find
$\delta>0,$ a sequence $(f_n)\in B_{\Cal F}$ and a sequence
$k_n\to \infty$ such that $k_n>2k_{n-1}$ and $\|T\hat
f_n(-k_n)e^{-k_n\theta}\|_{Y_{\theta}} \ge \delta $ for all
$n.$ Now, given $n$ and any $\epsilon>0$, we can use (a) to
find integers $m$ and $p$ such that $m>p\ge n$ and $\|T(\hat
f_m(-k_m)-\hat f_p(-k_p))\|_{Y_0} \le \epsilon.$ However
$\|T\hat f_m(-k_m)e^{-k_m}\|_{Y_1}\le 1$ and $\|T\hat
f_p(-k_p)e^{-k_m}\|_{Y_1} \le e^{k_p-k_m}\le 1$.  Hence,
again by $(2)$, $$ \|T(\hat f_m(-k_m)-\hat
f_p(-k_p))e^{-k_m\theta}\|_{Y_{\theta}}\le
C\epsilon^{1-\theta},$$ where $C$ depends only on $\theta $.
It follows that $$ \|T\hat
f_m(-k_m)e^{-k_m\theta}\|_{Y_{\theta}} \le
C(\epsilon^{1-\theta} + e^{(k_p-k_m)\theta}).$$ Hence $$
\delta \le C(\epsilon^{1-\theta} + e^{-k_n\theta})$$ and
this is a contradiction since $\epsilon>0$ and $k_n$ are
arbitrary.\bull\enddemo

\proclaim{Lemma 4}For $0<\theta<1$ and each fixed $N\in\bold
N$ the set $\{\Cal S_NTf(e^{\theta}):\ f\in B_{\Cal F}\}$ is
relatively compact in $Y_{\theta}.$\endproclaim

\demo{Proof}Suppose $f_n\in B_{\Cal F};$ then by Lemma 3(a)
we can pass to a subsequence $(g_n)$ such that for $|k|\le
2N$ we have $$ \|T\hat g_n(k) -T\hat g_{n+1}(k)\|_{Y_0} \le
2^{-n}.$$ Thus $$ \|\Cal S_NTg_n(z)-\Cal
S_NTg_{n+1}(z)\|_{Y_0}\le (4N+1)2^{-n}$$ for $|z|=1.$ Also,
for $|z|=e$, we have $$ \|\Cal S_NTg_n(z) -\Cal
S_NTg_{n+1}(z)\|_{Y_1} \le C_1$$ for some suitable constant
$C_1.$ Thus, by $(1)$, $$ \|\Cal S_NTg_n(e^{\theta}) -\Cal
S_NTg_{n+1}(e^{\theta})\|_{Y_{\theta}} \le C
C_1^{\theta}(4N+1)^{1-\theta}2^{-n(1-\theta)}$$ and so $\Cal
S_NTg_n(e^{\theta})$ is convergent.\bull\enddemo

Let $\Cal E$ be a subset of $\Cal F$.  We shall say that
$\Cal E$ is {\it effective} if it is bounded in $\Cal F$ and
if for some absolute constant $\lambda$ and every $f \in
\Cal E$ and every $n\in\bold N$ we have $f - \Cal S_nf \in
\lambda \Cal E$.  For each $\theta \in (0,1)$ let $\Cal
E_{\theta} = \{f(\theta):f\in \Cal E \}$.  We shall say that
$\Cal E$ is $\theta-${\it effective} if it is effective and
if $\Cal E_{\theta} \cap \gamma B_{X_{\theta}}$ is norm
dense in $\gamma B_{X_{\theta}}$ for some positive constant
$\gamma$ (which may depend on $\theta$).  Of course $B_{\Cal
F}$ is $\theta-$effective, but there are also clearly
smaller sets with the same property, for example the set of
those $f$ in $B_{\Cal F}$ with finitely many non zero
coefficients $\hat f(k)$.  (Cf.  [8] Section 9.2 and 29.2.)

\proclaim{Lemma 5}Let $\Cal E$ be an effective subset of
$\Cal F$ and let $\theta \in (0,1)$.  Then the following
conditions are equivalent:

(a) $T(\Cal E_{\theta})$ is a relatively compact subset of
$Y_{\theta}$.

(b) Every sequence $(f_n)$ in $\Cal E$ satisfies
$$\lim_{n\to\infty}\|Tf_n(e^{\theta})-\Cal
S_nTf_n(e^{\theta})\|_{Y_{\theta}}=0.  $$

If $\Cal E$ is $\theta-$effective then the preceding two
conditions are also equivalent to

(c) $T:X_{\theta}\to Y_{\theta}$ is compact.

\endproclaim

\demo{Proof}First suppose that (a) holds.  If $(f_n)$ is a
sequence in $\Cal E$ we observe that for a suitable constant
$C$ depending on $\theta$ we have $$ \|\Cal R(f_n -\Cal
S_nf_n)(e^{\theta})\|_{X_0+X_1} \le
Ce^{-n(1-\theta)}\max_{|z|=e}\|f_n(z)\|_{X_0+X_1}.$$
Combining this with a similar estimate for $\Cal
R_-(f_n-\Cal S_nf)(e^{\theta})$ we have $$ \|
f_n(e^{\theta})- \Cal S_nf_n(e^{\theta})\|_{X_0+X_1} \le
C(e^{-n(1-\theta)}+e^{-n\theta}).$$

Hence $$ \lim_{n\to\infty} \|Tf_n(e^{\theta})-\Cal S_n
Tf_n(e^{\theta})\|_{Y_0+Y_1}=0.$$ Using the fact that
$f_n(e^{\theta})-\Cal S_nf_n(e^{\theta}\} \in \lambda \Cal
E_{\theta}$ for each $n$ and condition (a) we deduce that we
also have convergence in $Y_{\theta}$, establishing (b).

Conversely, notice that if (b) holds then $\lim_{n\to\infty}
\|Tf(e^{\theta})-\Cal S_nTf(e^{\theta})\|_{Y_{\theta}}=0$
uniformly for $f\in \Cal E.$ It then follows from Lemma 4
that the set $\{Tf(e^{\theta}):\ f\in \Cal E\}$ is
relatively compact in $Y_{\theta}$ and so (a) holds.

Obviously (c) implies (a).  The reverse implication is also
trivial whenever $\Cal E$ is
$\theta-$effective.\bull\enddemo

\vskip2truecm

\subheading{3.  The main results}

The following theorem will imply the compactness result when
the domain space is $\langle X_0,X_1 \rangle_{\theta}$ or
$\langle X_0,X_1,\theta \rangle$ or when $\bold X$ is a
couple of lattices.

\proclaim{Theorem 6} Suppose that $\bold X$ and $\bold Y$
are Banach couples and that $T:\bold X\to\bold Y$ is such
that $T:X_0\to Y_0$ is compact.  Let $\Cal E$ be the subset
of $B_{\Cal F(\bold X)}$ consisting of those elements $f$
for which the series $\sum_{k\in\bold Z}e^{jk}\hat f(k)$
converges unconditionally in $X_j$ for $j=0,1$ and $\|
\sum_{k\in\bold Z}\lambda_k e^{jk}\hat f(k)\|_{X_0} < 1$ for
every sequence of complex scalars $(\lambda_k)$ with
$|\lambda_k| \le 1$ for all $k$.  Then $T(\Cal E_{\theta})$
is relatively compact in $Y_{\theta}$ for every
$0<\theta<1$.  \endproclaim

\demo{Proof}We may suppose that $\|T\|_{X_j \to Y_j}\le 1$
for $j=0,1$.  Consider an arbitrary sequence $(f_n)$ in
$\Cal E$ such that $\hat f_n(k)=0$ for $|k|\le n.$ Fix any
$0<\theta<1.$ Clearly $\Cal E$ is effective, so by Lemma 5
it will suffice to show that
$\lim_{n\to\infty}\|Tf_n(e^{\theta})\|_{Y_{\theta}}=0.$

For any $N\in\bold N$ let us pick a subset $A_n(N)$ of
$\bold Z$ so that $|A_n(N)|=N$ and $\|T\hat f_n(k)\|_{Y_0}
\le \|T\hat f_n(l)\|_{Y_0}$ whenever $k\notin A_n(N)$ and
$l\in A_n(N).$ Appealing to Lemma 3(d) we see that for any
fixed $N$ we must have $$ \lim_{n\to\infty}\|\sum_{k\in
A_n(N)}T\hat f_n(k)e^{k\theta}\|_{Y_{\theta}}=0.$$ It is
therefore possible to pick a non decreasing sequence of
integers $N_n$ with $\lim N_n=\infty$ so that $$
\lim_{n\to\infty}\|\sum_{k\in A_n(N_n)}T\hat
f_n(k)e^{k\theta}\|_{Y_{\theta}}=0.$$ We define $g_n(z)=
\sum_{k\notin A_n(N_n)}\hat f_n(k)z^k.$ Then it is easy to
check that $g_n\in B_{\Cal F}$.  Further, if
$b_n=\sup_{k\in\bold Z}\|T\hat g_n(k)\|_{Y_0}$, then
$\lim_{n\to\infty}b_n=0$ by Lemma 3(c).  It remains only to
show that
$\lim_{n\to\infty}\|Tg_n(e^{\theta})\|_{Y_{\theta}}=0.$

To this end suppose $y^*\in B_{Y_0^*}.$ Then, $$ \align
\int_0^{2\pi} |\langle Tg_n(e^{it}),y^*\rangle|^2
\frac{dt}{2\pi} &= \sum_{k\in\bold Z} |\langle T\hat g_n(k),
y^*\rangle|^2\\ &\le b_n \sum_{k\in\bold Z} |\langle T\hat
f_n(k),y^*\rangle|\\ &= b_n \sup_{|\lambda_k|\le 1}|\langle
T(\sum_{k\in\bold Z} \lambda_k f_n(k)),y^*\rangle|\\ &\le
b_n.  \endalign $$

Now by Lemma 1, $$ \lim_{n\to\infty}\int_0^{2\pi}
\|Tg_n(e^{it})\|_{Y_0}^2 \frac{dt}{2\pi}=0.$$ Finally we can
appeal to Lemma 2 (i) to obtain
$\lim_{n\to\infty}\|Tg_n(e^{\theta})\|_{Y_{\theta}}=0.$ This
completes the proof.\bull\enddemo

\proclaim{Corollary 7}For $\bold X, \bold Y$ and $T$ as
above,

(a) $T:\langle X_0,X_1 \rangle_{\theta} \to Y_{\theta}$ is
compact.

(b) $T:\langle X_0,X_1,\theta \rangle \to Y_{\theta}$ is
compact.

Furthermore if $\bold X$ is a couple of complexified Banach
lattices of measurable functions on some measure space then

(c) $T:X_{\theta} \to Y_{\theta}$ is compact.  \endproclaim
\demo{Proof}As pointed out [20] and in [16], $\langle
X_0,X_1 \rangle_{\theta}$ is contained in $X_{\theta}$.
More specifically, we observe that for each series
$x=\sum_{k\in\bold Z}x_k$ arising in the definition of an
element $x \in \langle X_0,X_1 \rangle_{\theta}$ we have
$\lim_{N \to \infty}\sup \|\sum_{|k|\ge N}\lambda_k
e^{(j-\theta)k}x_k\|_{X_j} =0$ for $j=0,1$ where the
supremum is over all choices of $\lambda_k$ with moduli $\le
1$.  Thus the function $f(z)=\sum_{k\in \bold Z}e^{-\theta
k}x_k z^k$ is $X_j- $continuous on $|z|=e^j$ and so it is an
element of $\Cal F(\bold X)$.  Consequently $\Cal
E_{\theta}$ is the open unit ball of $\langle X_0,X_1
\rangle_{\theta}$.  This immediately implies (a).

For (b) let $x$ be an arbitrary element in the open unit
ball of $\langle X_0,X_1,\theta \rangle$.  Then there exists
a representation $x=\sum_{k\in\bold Z}x_k$ for which the
elements $u_N = x=\sum_{|k|\le N}x_k$ are all in $\Cal
E_{\theta}$.  So by Theorem 6 there exists a subsequence of
$(Tu_N)$ which converges in the norm of $Y_{\theta}$ to some
element in the closure of $T(\Cal E_{\theta})$.  Since $u_N
\to x$ in $X_0+X_1$ this element must be $Tx$ and we deduce
that (b) holds.

If $\bold X$ is a couple of complexified Banach lattices
then $\langle X_0,X_1 \rangle_{\theta} = X_{\theta}$, as
follows from [8], section 13.6 (ii) p.125 and [22] Lemma
8.2.1 p. 453.  This of course establishes (c).
\bull\enddemo

Before proving the next theorem we will need a preliminary
lemma.

\proclaim{Lemma 8}Let $X$ be a UMD-space and let $V:X\to Y$
be a compact linear operator for some Banach space $Y$.
Then there exists a function $\eta:[0,\infty)\to [0,\infty)$
with $\lim_{\delta\to 0}\eta(\delta)=\eta(0)=0,$ such that
$\|\Cal R_-V\phi\|_{L_2(\bold T,Y)} \le
\eta(\|V\phi\|_{L_2(\bold T,Y)})$ for all $\phi\in L_2(\bold
T,X)$ with $\|\phi\|_{L_2(\bold T,X)}\le 2$.  \endproclaim

\demo{Proof} If the result is false then there is a sequence
$(\phi_n)$ for which $\|\phi_n\|_{L_2(\bold T,X)}\le 2,$
$\lim\|V\phi_n\|_{L_2(\bold T,Y)}=0$ but such that for some
$\epsilon>0$ we have $ \|\Cal R_-V\phi_n\|_{L_2(\bold T,Y)}
\ge \epsilon.$ However for all $y^*\in Y^*$ we have $$
\int_0^{2\pi}|\langle \Cal
R_-V\phi_n,y^*\rangle|^2\frac{dt}{2\pi} \le
\int_0^{2\pi}|\langle V\phi_n,y^*\rangle|^2\frac{dt}{2\pi}$$
and $\|\Cal R_-\phi_n\|_{L_2(\bold T,X)}$ is bounded by the
UMD-property of $X.$ Hence by Lemma 1 we obtain a
contradiction.\bull\enddemo

\proclaim{Theorem 9}Suppose that $\bold X=(X_0,X_1)$ and
$Y=(Y_0,Y_1)$ are Banach couples and $X_0$ is a UMD-space.
Let $T:\bold X\to\bold Y$ be such that $T:X_0\to Y_0$ is
compact.  Then $T:X_{\theta}\to Y_{\theta}$ is compact for
every $0<\theta<1$.\endproclaim

\demo{Proof}Using Lemma 2 (iii) we see that we may assume
without loss of generality that both of the couples $\bold
X$ and $\bold Y$ are regular.  This ensures that the dual
spaces also form Banach couples.  In particular we will make
use of the fact that $(Y_0+Y_1)^* = Y_0^*\cap Y_1^*$ (cf.
[3] p. 32) and so this space clearly separates points of
$Y_0$.

As in the proofs of preceding theorems, it will suffice to
consider a sequence $f_n\in B_{\Cal F}$ satisfying the
conditions $\hat f_n(k)=0$ for $|k|\le n$ and show that
$\lim_{n\to\infty}\|Tf_n(e^{\theta})\|_{Y_{\theta}}=0.$ We
may of course suppose as before that $\|T\|_{\bold X \to
\bold Y} \le 1.$

We first consider $\Cal Rf_n$.  We note that for $|z|=1$ we
have an estimate $ \|\Cal Rf_n(z)\|_{X_1}\le Ce^{-n}$ and,
by the UMD-property of $X_0$, the sequence $\Cal Rf_n$ is
bounded in $L_2(\bold T,X_0).$ For each $y^*\in Y_1^*\cap
Y_0^*$ we see that $\langle \Cal RTf_n,y^*\rangle$ is
uniformly convergent to $0$.  So we can apply Lemma 0 to
deduce that $$ \lim_{n\to\infty} \|\Cal RTf_n\|_{L_2(\bold
T,Y_0)} =0.  \leqno{(6)} $$

Let us fix $\epsilon>0.$ Since $T:X_0 \to Y_0$ is compact
and $X_0\cap X_1$ is dense in $X_0$ we can find a finite set
$\{x_1,x_2,\ldots,x_N\}$ in $B_{X_0}\cap X_1$ so that if
$\|x\|_{X_0}\le 1$ then there exists $1\le j\le N$ with
$\|Tx-Tx_j\|_{Y_0}\le \epsilon.$ Thus for each $n$ we can
find a measurable function $H_n:\bold T\to
\{x_1,\ldots,x_n\}$ so that
$\|Tf_n(e^{it})-TH_n(e^{it})\|_{Y_0}\le \epsilon$ for all
$t.$ By convolving with a suitable kernel we can obtain a
$C^{\infty}-$function $h_n:\bold T\to F$ (where $F$ is the
linear span of $\{x_1,\ldots,x_n\}$) so that
$\|Tf_n(e^{it})-Th_n(e^{it})\|_{Y_0}\le 2\epsilon$ and
$\|h_n(e^{it})\|_{X_0}\le 1$ for all $t$.  Let us expand
$h_n$ in its Fourier series $$h_n(e^{it}) =\sum_{n\in\bold
Z}\hat h_n(k)e^{ikt}.$$ We will define $$ g_n(z) =
\sum_{k\le-(n+1)}\hat h_n(k) z^k$$ for $|z|\ge 1.$ This
defines an $F$-valued function which is analytic for $|z|>1$
and continuous for $|z|\ge 1$ since $h_n$ is $C^{\infty}.$

Now $\Cal R_-f_n-g_n= z^{-n}\Cal R_-(z^nf_n-z^ng_n) =
z^{-n}\Cal R_-(z^nf_n-z^nh_n)$.  Also clearly the functions
$\phi_n = z^nf_n-z^nh_n$ satisfy $\| \phi _{n}(e^{it})\|
_{X_{0}}\le 2$ and so $\| \phi _{n}\| _{L^2(\bold T,X_0)}\le
2$.  Thus we can apply Lemma 8 to obtain that $\| {\Cal
R}_{-}Tf_{n}-Tg_{n}\| _{L^2(\bold T,Y_0)} = \| {\Cal
R}_{-}T\phi _{n}\| _{L^2(\bold T,Y_0)} \le \eta (2\epsilon
)$ for some function $\eta $ which depends only on $T$ and
satisfies $\lim_{\delta \to 0}\eta (\delta ) = 0$.
Combining this with $(6)$ we obtain that
$\limsup_{n\to\infty}\|Tf_n-Tg_n\|_{L_2(\bold T,Y_0)} \le
\eta (2\epsilon )$.

Now consider $(g_n)$ on the circle $|z|=e.$ For a suitable
constant $C_1$ we have $$\|g_n(z)\|_{X_0} \le
C_1e^{-n}\max_{|\zeta|=1}\|h_n(\zeta)\|_{X_0}\le
C_1e^{-n}.$$ Since $g_n$ is $F$-valued there is a constant
$C_2$, depending only on $F$ and thus on $\epsilon$, such
that $\|x\|_{X_1}\le C_2\|x\|_{X_0}$ for all $x\in F.$ Thus
we have $\lim_{n\to\infty}\max_{|z|=e}\|g_n(z)\|_{X_1}=0.$
>From this we conclude that
$\limsup_{n\to\infty}\max_{|z|=e}\|Tf_n(z)-Tg_n(z)\|_{Y_1}
\le 1.$

Now we can deduce, using Lemma 2 (i), that
$$\limsup_{n\to\infty}\|Tf_n(e^{\theta})
-Tg_n(e^{\theta})\|_{Y_{\theta}} \le C_3(\eta (2\epsilon
))^{\theta}$$ for a constant $C_3$ which depends only on
$\theta$.

However we can also estimate $\|g_n(e^{\theta})\|_{X_0} \le
C_4e^{-n\theta}$ and again using the fact that all $g_n$
have range in $F$ we have $\lim
\|g_n(e^{\theta})\|_{X_{\theta}}=0.$ Thus we are left with
the estimate $$
\limsup_{n\to\infty}\|Tf_n(e^{\theta})\|_{Y_{\theta}}\le
C_3(\eta (2\epsilon ))^{\theta}$$ Since $\epsilon>0$ is
arbitrary this completes the proof.\bull\enddemo

\demo{Remark}See the introduction for a discussion of the
class of UMD-spaces.  \enddemo

\proclaim{Theorem 10}Let $\bold X$ be a Banach couple such
that $X_0$ is reflexive and is given by
$X_0=[W,X_1]_{\alpha}$ for some $0<\alpha<1$ and some Banach
space $W$ which forms a Banach couple with $X_1$.  Suppose
$T:\bold X\to \bold Y$ is such that $T:X_0\to Y_0$ is
compact.  Then $T:X_{\theta}\to Y_{\theta}$ is compact for
$0<\theta<1.$\endproclaim

\demo{Proof}Letting the notation $X^{\circ }$ now mean the
closure of $W\cap X_{1}$ in $X$, we observe that $[W^{\circ
},X^{\circ }_{1}]_{\delta}= [W,X_{1}]_{\delta }$ for all
$\delta \in (0,1)$ (cf.  Lemma 2 (iii)).  Consequently we
may assume without loss of generality that $W\cap X_{1}$ is
dense in $W$ and also in $X_{1}$.

For each $\theta \in (0,1)$ we have $$X_{\theta
}=[[W,X_{1}]_{\alpha },X_{1}]_{\theta }= [W,X_{1}]_{\delta
}$$ where $\delta =(1-\theta )\alpha +\theta $. (Cf.  Lemma
2 (iv).)

Let $\Cal E$ be the set consisting of all functions in
$B_{\Cal F(\bold X)}$ which can be extended to functions $f$
on the closed annulus $\{z:e^{-\beta}\le |z|\le e\}$ where
$\beta=\alpha/(1-\alpha)$ in such a way that $f$ is analytic
into $W+X_1$ on the open annulus, $W+X_1-$continuous on the
closed annulus, $W-$continuous on $|z|=e^{-\beta}$ and
$\max_{|z|=e^{-\beta}}\|f(z)\|_W \le 1.$ Clearly $\Cal E$ is
effective.  Furthermore it is also $\theta-$effective for
every $\theta \in (0,1)$.  This can be shown readily using
the above reiteration formula together with the observation
[11] that the complex interpolation method yields the same
spaces on annuli of different dimensions, even if they are
not conformally equivalent.  (The spaces defined using the
strips $\{z:0\le \Re{z}\le 1\}$ and $\{z:-\beta \le
\Re{z}\le 1\}$ are obviously identical.  Now simply
``periodize" the functions on both of these strips with
period $2\pi i$ as in [11].)

We will apply Lemma 5. Consider a sequence $f_n\in\Cal E.$
Let $g_n=f_n-\Cal S_nf_n.$ Suppose $x^*\in W^*\cap X_1^*.$
Then $$ \align \int_0^{2\pi}|\langle
g_n(e^{it}),x^*\rangle|^2 \frac{dt}{2\pi} &= \sum_{k\in\bold
Z}|\langle \hat g_n(k),x^*\rangle|^2\\ &\le \sum_{k\le
-n}|\langle \hat f_n(k),x^*\rangle|^2 +\sum_{k\ge n}|\langle
\hat f_n(k),x^*\rangle|^2\\ &\le C(e^{-2\beta
n}\|x^*\|_{W^*}^2 + e^{-2n}\|x^*\|_{X_1^*}^2).  \endalign $$
By our density assumption $W^* \cap X^*_1 = (W+X_1)^*$ so
this space separates points of $X_0 \subset W+X_1$.  It
follows that the set $U$ of $x^*\in X_0^*$ such that $$
\lim_{n\to\infty} \int_0^{2\pi}|\langle
g_n(e^{it}),x^*\rangle|^2 \frac{dt}{2\pi} =0$$ is a closed
weak$^*$ dense subspace of $X_0^*$.  Since $X_0$ is
reflexive $U=X_0^*$ and so $T^*(Y_0^*)\subset U$.  By Lemma
1 we obtain that $\lim_{n\to\infty}\|Tg_n\|_{L_2(\bold
T,Y_0)}=0$ and then an application of Lemma 2 (i) gives that
$\lim_{n\to\infty}\|Tg_n(e^{\theta})\|_{Y_{\theta}}=0.\bull$
\enddemo

\demo{Remark}The reader may care to note that if the
preceding theorem can be proved without the requirement that
$X_0$ is reflexive then Question 1 is completely answered
for the complex method, by using the reduction of this
problem given in [12] p. 339 to the case where ${\bold
X}=(l_1(FL_1),l_1(FL_1(e^{\nu}))$ and ${\bold
Y}=(l_{\infty}(FL_{\infty}),l_{\infty}(FL_{\infty}(e^{\nu}))$.
In this case we can of course take $W = l_1(FL_1(e^{-\beta
\nu})$ for $\beta $ as above.  \enddemo

Here is a sort of ``dual" result to Theorem 10.  Note that
it does not require any reflexivity conditions.  But
unfortunately it is still not sufficient to give a complete
answer to Question 1 (cf. the preceding remark) since
$[l_{\infty}(FL_{\infty}(e^{-\beta
\nu}),l_{\infty}(FL_{\infty}(e^{\nu}) ]_{\alpha }$ is
strictly contained in $l_{\infty}(FL_{\infty})$.

\proclaim{Theorem 11}Suppose $T:{\bold X}\to {\bold Y}$
where $T:X_{0}\to Y_{0}$ is compact.  Suppose that for some
Banach space $Z$, $(Z,Y_{1})$ forms a Banach couple and
$Y_{0}=[Z,Y_{1}]_{\alpha }$ for some $\alpha \in (0,1)$.
Then $T:X_{\theta }\to Y_{\theta }$ is compact for each
$\theta \in (0,1)$.  \endproclaim

\demo{Proof}We begin by showing that we can reduce the proof
to the case where a number of density conditions are
satisfied.  First, using Lemma 2 (iii) and rather similar
reasoning to before, we can suppose without loss of
generality that $\bold X$ is regular, and similarly, that
$Y_{0}\cap Y_{1}$ is dense in $Y_{1}$.  (The hypotheses
already ensure that $Y_{0}\cap Y_{1}$ is dense in $Y_{0}$.)
In fact we can furthermore suppose that $Z\cap Y_{1}$ is
dense in $Y_{1}$, since if that were not so we could replace
the couples $\bold X=(X_{0},X_{1})$ and $\bold
Y=(Y_{0},Y_{1})$ by $(X_{0},X_{\sigma })$ and
$(Y_{0},Y_{\sigma })$ for some number $\sigma \in (\theta
,1)$.  By several applications of Lemma 2 (iv) these latter
couples also satisfy all the other required hypotheses of
the theorem and we will be able to deduce the original
desired conclusion for $T:X_{\theta }\to Y_{\theta }$ since
$X_{\theta }=[X_{0},X_{\sigma }]_{\theta /\sigma }$ and
$Y_{\theta }=[Y_{0},Y_{\sigma }]_{\theta /\sigma }$.
Finally, given that all the above density conditions hold,
we can now, if necessary, replace $Z$ by $Z^{\circ}$, the
closure of $Z\cap Y_1$ in $Z$ without changing any of the
other spaces.  Also of course $Z^{\circ}\cap Y_1$ is dense
in $Y_1$.  In other words, we can also assume that $Z\cap
Y_1$ is dense in $Z$.

Let $T^*:(Y^*_0+Y^*_1) \to (X^*_0+X^*_1)$ be the adjoint of
$T:X_0 \cap X_1 \to Y_0 \cap Y_1$.  Clearly $T^*$ maps
$Y^*_{1}$ to $X^*_{1}$ boundedly and $Y^*_{0}$ to $X^*_{0}$
compactly.  This means that $T^*:[Z^*,Y^*_{1}]_{\alpha }\to
X^*_{0}$ is compact, since by Calder\'on's duality theorem
$Y^*_{0}=[Z^*,Y^*_{1}]^{\alpha }$ and $[Z^*,Y^*_{1}]_{\alpha
}$ is a closed subspace of $[Z^*,Y^*_1]^{\alpha }$.  (See
[2].)  Thus the operator $T^*$ satisfies all the hypotheses
of Theorem 10, ($T^*$ replaces $T$, $Z^*$ plays the role of
$W$, and instead of the original couples ${\bold X}$ and
${\bold Y}$ we have ${\bold Z}^*=([Z^*,Y^*_{1}]_{\alpha
},Y^*_{1})$ and ${\bold X}^*=(X^*_{0},X^*_{1})$
respectively) except that $[Z^*,Y^*_{1}]_{\alpha }$ is not
necessarily reflexive.

We now define ${\cal E}$ exactly analogously to the
definition in the proof of Theorem 10, i.e. it is the subset
of $B_{{\Cal F}({\bold Z}^{*_{)}}}$ of functions which are
extendable to $Z^*+Y^*_{1}-$valued continuous functions on
the annulus $\{z:e^{-\beta } \le \mid z\mid \le e\}$ which
are analytic in the interior of the annulus and are
continuous into $Z^*$, respectively $Y^*_1$ on the inner,
respectively outer components of the boundary.  Again we
consider the sequence $g_{n}=f_{n}-{\cal S}_{n}f_{n}$, where
$f_{n}$ is an arbitrary sequence in ${\Cal E}$.

This time we let $U$ be the set of all $y \in Y_0$ such that
$$ \lim_{n\to\infty} \int_0^{2\pi}|\langle
y,g_n(e^{it})\rangle|^2 \frac{dt}{2\pi} =0.$$ Using
estimates similar to those in the proof of Theorem 10 we
obtain that $Z\cap Y_1\subset U$.  Since $U$ must be closed
in $Y_0$ it follows that $U=Y_{0}$.  Consequently, $$
\lim_{n\to\infty} \int_0^{2\pi}|\langle
x,T^*g_n(e^{it})\rangle|^2 \frac{dt}{2\pi} =0$$ for all $x
\in X_0$.  We can now apply Lemma 1 to
$T^*:[Z^*,Y^*_1]_{\alpha } \to X^*_0$ to obtain that
$\lim_{n\to\infty}\|T^*g_n\|_{L_2(\bold T,X^*_0)}=0$.  Then
Lemma 2 gives that
$\lim_{n\to\infty}\|T^*g_n(e^{\theta})\|_{[X^*_0,X^*_1]_{\theta}}=0.$
By Lemma 5 we deduce that
$T^*:[[Z^*,Y^*_1]_{\alpha},Y^*_1]_{\theta} \to
[X^*_0,X^*_1]_{\theta}$ is compact.

As already remarked above, $[Z^*,Y^*_1]_{\alpha }$ is a
closed subspace of $Y^*_0$.  Furthermore,
$[Z^*,Y^*_1]_{\alpha }$ contains $Z^*\cap Y^*_1$ densely and
so obviously it is also the closure in $Y^*_0$ of the larger
space $[Z^*,Y^*_1]_{\alpha } \cap Y^*_1$.  So Lemma 2 (iii)
yields that $[[Z^*,Y^*_1]_{\alpha},Y^*_1]_{\theta}=
[Y^*_0,Y^*_1]_{\theta}$.

Let $z^*$ be an arbitrary element of the open unit ball of
$[Y^*_0,Y^*_1]^{\theta }$.  Thus $z^{*}= h^\prime (\theta )$
where $h$ is an element of the unit ball of the space
$\overline {{\Cal F}}(Y^{*}_{0},Y^{*}_{1})$ (of analytic
functions on the {\it unit strip} as defined in [8]).  If we
set $h_{n}(z) = ne^{(z^2-1)/n}(h(z+1/n)-h(z))$, and
$y^{*}_{n}=h_{n}(\theta )$ then it is easy to see that
$(y^{*}_{n})$ is a sequence in the unit ball of
$[Y^{*}_{0},Y^{*}_{1}]_{\theta }$ which converges to $z^{*}$
in $Y^{*}_{0}+Y^{*}_{1}$.  (Cf.  [11] p. 1006.)  In view of
the compactness of $T^{*}:[Y^{*}_{0},Y^{*}_{1}]_{\theta }\to
[X^{*}_{0},X^{*}_{1}]_{\theta }$ we can suppose that (some
subsequence of) the sequence $(T^{*}y^{*}_{n})$ is Cauchy in
$[X^{*}_{0},X^{*}_{1}]_{\theta }$.  Thus its limit in
$[X^{*}_{0},X^{*}_{1}]_{\theta }$ is also its limit in
$X^{*}_{0}+X^{*}_{1}$, namely $T^{*}z^{*}$.  This shows that
$T^{*}$ maps the unit ball of $[Y^{*}_{0},Y^{*}_{1}]^{\theta
}$ into a relatively compact subset of
$[X^{*}_{0},X^{*}_{1}]_{\theta }\subset
[X^{*}_{0},X^{*}_{1}]^{\theta }$.  Consequently
$T^{*}:[Y^{*}_{0},Y^{*}_{1}]^{\theta }\to
[X^{*}_{0},X^{*}_{1}]^{\theta }$ is compact.  This, together
with Calder\'on's duality theorem and the classical Schauder
theorem, completes the proof.  \bull \enddemo

\subheading{Acknowledgements} We are grateful to Natan
Krugljak, Mieczys\l aw Masty\l o and Rich\-ard Rochberg for
many very helpful discussions related to the problem of
interpolation of compact operators by the complex method.

\vskip2truecm

\subheading{References}

\item{1.}C.  Bennett and R. Sharpley, {\it Interpolation of
Operators,} Academic Press, New York 1988.

\item{2.}J.  Bergh, On the relation between the two complex
methods of interpolation, Indiana Univ.  Math.  J. 28
(1979), 775-778.

\item{3.}J.  Bergh and J. L\"ofstr\"om, {\it Interpolation
spaces.  An Introduction,} Springer, Berlin 1976.

\item{4.}J.  Bourgain, Some remarks on Banach spaces in
which martingale difference sequences are unconditional,
Ark.  Mat. 21 (1983) 163-168.

\item{5.}Y.  Brudnyi and N. Krugljak, {\it Interpolation
functors and interpolation spaces, Volume 1,} North Holland,
Amsterdam 1991.

\item{6.}D.L.  Burkholder, A geometrical characterization of
Banach spaces in which martingale difference sequences are
unconditional, Ann.  Prob. 9 (1981) 997-1011.

\item{7.}D.L.  Burkholder, A geometric condition that
implies the existence of certain singular integrals of
Banach space-valued functions, 270-286 in {\it Conference on
harmonic analysis in honor of A. Zygmund,} (W.  Beckner,
A.P.  Calder\'on, R. Fefferman, P.W.  Jones, editors)
Wadsaworth, Belmont, California, 1983.

\item{8.}A.  P. Calder\'on, Intermediate spaces and
interpolation, the complex method, Studia Math. 24 (1964)
113-190.

\item{9.}F.  Cobos, T. K\"uhn and T. Schonbek, One-sided
compactness results for Aronszajn-Gagliardo functors.  J.
Functional Analysis 106 (1992), 274-313.

\item{10.}F.  Cobos and J. Peetre, Interpolation of
compactness using Aronszajn-Gagliardo functors.  Israel J.
Math. 68 (1989), 220-240.

\item{11.}M.  Cwikel, Complex interpolation, a discrete
definition and reiteration, Indiana Univ.  Math.  J. 27,
1978, 1005-1009.

\item{12.}M.  Cwikel, Real and complex interpolation and
extrapolation of compact operators.  Duke Math.  J. 65
(1992) 333-343.

\item{13.}J.  Gustavsson and J. Peetre, Interpolation of
Orlicz spaces, Studia Math. 60 (1977), 33-59.

\item{14.}J.  A. Gutierrez, {\it On the boundedness of the
Banach space-valued Hilbert transform,} Ph.D. thesis,
University of Texas, Austin 1982.

\item{15.}K.  Hayakawa, Interpolation by the real method
preserves compactness of operators, J. Math.  Soc.  Japan 21
(1969), 189-199.

\item{16.}S.  Janson, Minimal and maximal methods of
interpolation, J. Functional Analysis 44 (1981), 50-73.

\item{17.}M.  A. Krasnosel'ski\v{i}, On a theorem of M.
Riesz, Soviet Math.  Dokl. 1 (1960), 229-231.

\item{18.}J.  L. Lions and J. Peetre, Sur une classe
d'espaces d'interpolation.  Inst.  Hautes Etudes Sci.  Publ.
Math. 19 (1964), 5-68.

\item{19.}M.  Masty\l o, On interpolation of compact
operators.  Preprint.

\item{20.}J.  Peetre, Sur l'utilization des suites
inconditionellement sommables dans la th\'eorie des espaces
d'interpolation.  Rend.  Sem.  Mat. Univ.  Padova 46 (1971),
173-190.

\item{21.}A.  Persson, Compact linear mappings between
interpolation spaces.  Ark.  Mat. 5 (1964), 215-219.

\item{22.}V.I.  Ovchinnikov, The method of orbits in
interpolation theory, Mathematical Reports, Vol 1, Part 2,
Harwood Academic Publishers 1984, 349-516.

\bye